\documentclass[11pt]{article}
\usepackage[dvips]{graphics}
\usepackage{amsmath}
\usepackage{amssymb}
\usepackage{amsfonts}

\def\a{\alpha}
\def\b{\beta}

\def\ga{\Gamma}

\def\l{\lambda}

\def\o{\omega}
\def\om{\Omega}

\def\iy{\infty}
\def\im{{\rm Im\, }}

\def\BD{{\mathbb D}}

\newcommand{\sH}{{\mathcal H}}
\newcommand{\sK}{{\mathcal K}}
\newcommand{\sN}{{\mathcal N}}

\newcommand{\sF}{{\mathcal F}}
\newcommand{\sU}{{\mathcal U}}

\newcommand{\sG}{{\mathcal G}}
\newcommand{\sL}{{\mathcal L}}
\newcommand{\sE}{{\mathcal E}}
\newcommand{\sD}{{\mathcal D}}

\newcommand{\sY}{{\mathcal Y}}
\newcommand{\sX}{{\mathcal X}}

\newcommand{\sZ}{{\mathcal X}}

\newcommand{\ov}[1]{{\overline{#1}}}

\newcommand{\eH}{{\mathbf{H}}}
\newcommand{\eS}{{\mathbf{S}}}

\newcommand{\bpr}{{\noindent\textbf{Proof.}\ \ }}
\newcommand{\epr}{{\hfill $\Box$}}

\newtheorem{thm}{Theorem}[section]
\newtheorem{lem}[thm]{Lemma}
\newtheorem{cor}[thm]{Corollary}
\newtheorem{prop}[thm]{Proposition}

\title
{Relaxed commutant lifting: an  equivalent version and a new
application}
 \author
{A.E. Frazho, S. ter Horst and M.A. Kaashoek}

\date{}

\begin{document}
\maketitle \thispagestyle{empty}

\begin{abstract}
This paper presents a few additions to commutant lifting theory. An
operator interpolation problem is introduced and shown to be
equivalent to the relaxed commutant lifting problem. Using this
connection a description of all solutions of the former problem is
given. Also a new application, involving bounded operators induced
by $H^2$ operator-valued functions, is presented.
\end{abstract}

\setcounter{section}{-1}
\section{Introduction}
Let $\sU$ and $\sY$ be Hilbert spaces, and let $\sF$ be a subspace
of $\sU$. In this paper we consider the following problem. Given a
contraction
\begin{equation}\label{defom}
\o=\left[\begin{array}{c}\o_1\\
\o_2
\end{array}\right]:\sF\to \left[\begin{array}{c}\sY\\
\sU
\end{array}\right],
\end{equation}
find   a (all) contraction(s) $\ga$ from $\sU$ into $H^2(\sY)$
satisfying the equation
\begin{equation}\label{cond1}
E_\sY\o_1 +S_\sY\ga \o_2 =\ga|_\sF.
\end{equation}
Here and in the sequel we use the convention that for any Hilbert
space $\sN$ the symbol $S_\sN$ denotes the forward shift on the
Hardy space $H^2(\sN)$ and $E_\sN$ denotes the embedding of $\sN$
into $H^2(\sN)$ defined by $(E_\sN\, n)(\l)\equiv n$. Furthermore,
$\ga|_\sF$ stands for the restriction of $\ga$ to $\sF$ viewed as
an operator from $\sF$ into $H^2(\sY)$.

A contraction $\ga$ from $\sU$ into $H^2(\sY)$ satisfying equation
(\ref{cond1}) will be called a \emph{solution to the interpolation
problem defined by the contraction} $\o$ in (\ref{defom}).

We shall show that the problem stated above can be reformulated as
a relaxed commutant lifting problem. On the other hand, as we know
from \cite{FFK02a}, the relaxed commutant lifting problem can be
reduced to an interpolation problem defined by a special
contraction $\o$ of the form (\ref{defom}). Hence it follows that
the above problem and the relaxed commutant lifting problem are
equivalent in the sense that the one problem can be reduced to the
other and conversely.

To state the main results we need some additional notation.
Throughout $\sL(\sU,\sY)$ stands for the space of all (bounded
linear) operators from  $\sU$ into $\sY$. By $\eH^2(\sU,\sY)$ we
denote the space of all $\sL(\sU,\sY)$-valued functions that are
analytic on $\BD$ such that the Taylor coefficients
$H_0,H_1,H_2,\ldots$ of the function $H$ at zero satisfy the
constraint $\sum_{n=0}^\infty \|H_nu\|^2<\infty$ for each $u\in
\sU$. Given such a function $H$, the formula
\begin{equation}
\label{gaH} (\ga u)(\l)=H(\l)u, \quad \l\in \BD, \quad u\in \sU,
\end{equation}
defines an operator $\ga$ from $\sU$ into the Hardy space
$H^2(\sY)$, which we shall refer to as  the \emph{operator defined
by} $H$. Conversely, if $\ga$ is an operator from $\sU$ into
$H^2(\sY)$, then there is a unique $H\in \eH^2(\sU,\sY)$ such that
(\ref{gaH}) holds, and in this case we call $H$ the \emph{defining
function} of $\ga$.

Replacing $\ga$ in (\ref{cond1}) by its defining function we see
that our problem has the following alternative formulation: find
all $H\in \eH^2(\sU,\sY)$ satisfying
\begin{equation}\label{cond2}
\o_1 +\l H(\l) \o_2 =H(\l)|_\sF,\quad \l\in \BD,
\end{equation}
and such that the operator defined by $H$ is a contraction. In
this case we also say that  $H$ is   a \emph{solution to the
interpolation problem defined by the contraction} $\o$ in
(\ref{defom}).

The connection with  relaxed commutant lifting mentioned above
allows us to use Theorem 1.1 in \cite{FtHK06b} (cf., Theorem 0.1
in \cite{FtHK06a}) to prove the following theorem.

\begin{thm}
\label{thmmain1} Let $\o$ be a contraction as in $(\ref{defom})$.
Then $H\in \eH^2(\sU,\sY)$ is a solution to the interpolation
problem defined by the contraction $\o$ if and only if $H$ is given
by
\begin{equation}\label{sols2}
H(\l)=\Pi_{\sY}Z(\l)(I-\l\Pi_{\sU}Z(\l))^{-1}, \quad \l\in\BD,
\end{equation}
where $Z$ is an arbitrary Schur  class function from
$\eS(\sU,\sY\oplus\sU)$ satisfying the constraint
$Z(\l)|_{\sF}=\o$ for each $\l\in\BD$. Here $\Pi_{\sY}$ and
$\Pi_{\sU}$ are the orthogonal projections from the Hilbert space
direct sum $\sY\oplus\sU$ onto $\sY$ and $\sU$, respectively.
\end{thm}
Recall that for Hilbert spaces $\sH$ and $\sK$ the Schur class
$\eS(\sH,\sK)$ consists of all $\sL(\sH,\sK)$-valued functions
$F$, analytic on $\BD$, such that $\sup_{\l\in\BD}\|F(\l)\|\leq
1$.

In general, the map $Z\mapsto H$ defined by Theorem \ref{thmmain1}
is not one-to-one. In fact, using the connection with relaxed
commutant lifting and Theorem 1.2 in \cite{FtHK06b} we shall
derive the following result.

\begin{thm}
\label{thmmain2}Let $H\in \eH^2(\sU,\sY)$ be a solution to the
interpolation problem defined by the contraction $\o$ in
$(\ref{defom})$, and let $\ga $  from $\sU$ into $H^2(\sY)$ be the
operator defined by $H$. Then the set of all $Z\in
\eS(\sU,\sY\oplus\sU)$ satisfying   $Z(\l)|_{\sF}=\o$ for each
$\l\in\BD$ and such that $(\ref{sols2})$ holds is parameterized by
the set
\begin{equation}
\label{parsetC1} \big\{ C\in \eS(\sD_\ga, \sD_\ga)\mid
C(\l)D_\ga|_\sF=D_\ga\o_2\ \mbox{for each  $\l\in \BD$}\big\}.
\end{equation}
\end{thm}

The parameterization referred to in the preceding theorem can be
made more explicit. Indeed, let $H\in \eH^2(\sU,\sY)$ be a solution 
to the interpolation problem defined by the contraction $\o$ in (\ref{defom}), 
and let the operator $\ga$ defined by $H$ be a contraction. Then given 
$C\in\eS(\sD_\ga, \sD_\ga)$ we define
\begin{equation}
\label{defZC} Z_C(\l)=\left[\begin{array}{c} 2H(\l)(W(\l)+I)^{-1}\\
\noalign{\vskip6pt} \l^{-1}(W(\l)-I)
\end{array}\right], \quad\l\in\BD,
\end{equation}
where
\begin{eqnarray}\label{defW}
&&W(\l)=\ga^*(I+\l S_\sY^*)(I-\l S_\sY^*)^{-1}\ga+\nonumber\\
\noalign{\vskip6pt} &&\hspace{2.5cm}+ \,D_\ga(I+\l C(\l))(I-\l
C(\l))^{-1}D_\ga, \quad\l\in\BD.
\end{eqnarray}
We shall see that the map  $C\mapsto Z_C$ induces a one-to-one map
from the set (\ref{parsetC1}) onto the set of all $Z\in
\eS(\sU,\sY\oplus\sU)$ satisfying   $Z(\l)|_{\sF}=\o$ for each
$\l\in\BD$ and such that $(\ref{sols2})$ holds.

As a new application we shall use Theorem \ref{thmmain1} to prove
the following result.

\begin{thm}\label{thmrepr}
Let $H\in\eH^2(\sU,\sY)$, and let $\Theta\in \eS(\sE,\sU)$ be
inner such that  $\Theta(0)=0$. Put $\sH=H^2(\sU)\ominus\Theta
H^2(\sE)$. In order that the  map $f\mapsto Hf$ defines a
contraction from $\sH$ into $H^2(\sY)$ it is necessary and
sufficient that $H$ is given by
\begin{equation}
\label{functH}
H(\l)=\Pi_{\sY}Z(\l)(I-\Theta(\l)\Pi_{\sE}Z(\l))^{-1}, \quad
\l\in\BD,
\end{equation}
where $Z$ is an arbitrary Schur  class function from
$\eS(\sU,\sY\oplus\sE)$.
\end{thm}
For $\Theta(\l)=\l^N$ the matrix-valued version of the above
theorem can been found in \cite{ABL96, ABL97} and for
operator-valued functions in \cite{tH06a}. For the scalar case,
with  $\Theta(\l)=\l$, the result goes back to \cite{Sar89}, page
490.

\medskip
The paper consist of three sections (not counting the present
introduction). The first section has a preliminary character. Here
we recall how the relaxed commutant lifting problem can be reduced
to an interpolation problem of the type defined above. In the second
section we prove Theorems \ref{thmmain1} and \ref{thmmain2} using
relaxed commutant lifting. In the third section Theorem
\ref{thmrepr} is proved.

\section{Preliminaries about relaxed commutant lifting}
\setcounter{equation}{0}

This section has a preliminary character. We recall the relaxed
commutant lifting problem and how this problem can be reduced to
an interpolation problem defined by a contraction of the form
(\ref{defom}).

We begin with some terminology. A quintet
$\{A,T^\prime,U^\prime,R,Q\}$ consisting of five Hilbert space
operators is called a \emph{data set} if the operator $A$ is a
contraction mapping $\sH$ into ${\sH}'$, the operator $U^\prime$ on
$\sK^\prime$ is a minimal isometric lifting of the contraction
$T^\prime$ on $\sH'$, and $R$ and $Q$ are operators from $\sH_0$ to
$\sH$ satisfying the following constraints:
\begin{equation}\label{intertw}
T'AR=AQ  \quad \textup{and}\quad R^*R\leq Q^*Q.
\end{equation}
Without loss of generality we can and shall assume that $U^\prime$
is the Sz.-Nagy-Sch\"{a}ffer (minimal) isometric lifting of
$T^\prime$. The latter means (see \cite{Sz.-NF}) that
\begin{equation}\label{szns}
 U^\prime = \left[\begin{array}{cc}
   T^\prime & 0 \\
   E_{\sD_{T'}}D_{T^\prime} & S_{\sD_{T'}} \\
 \end{array}\right] \mbox{ on } \sK' =  \left[\begin{array}{c}
   \sH^\prime \\
   H^2(\sD_{T^\prime}) \\
 \end{array}\right].
\end{equation}

Given this data set the relaxed commutant lifting problem
(\emph{RCL problem}) is to find all contractions $B$ from $\sH$ to
$\sK'$ such that
\begin{equation}\label{rclt}
\Pi_{\sH'}B=A\quad\mbox{and}\quad U'BR=BQ.
\end{equation}
Here $\Pi_{\sH'}$ is the orthogonal projection from $\sK'$ onto
$\sH'$ viewed as an operator from $\sK'$ into $\sH'$. In this case
we refer to $B$ as a \emph{solution to the RCL problem} for the data
set $\{A,T^\prime,U^\prime,R,Q\}$.

Since $U^\prime$ is given by (\ref{szns}), we have $\sK'=\sH'\oplus
H^2(\sD_{T'})$, and an operator $B$ from $\sH$ into $\sH'\oplus
H^2(\sD_{T'})$ is a contraction satisfying the first identity in
(\ref{rclt}) if and only if $B$ can be represented in the form
\begin{equation}\label{defB}
B= \left[\begin{array}{c}
     A  \\
   \Gamma D_A \\
 \end{array}\right]: \sH \to \left[\begin{array}{c}
   \sH^\prime \\
   H^2(\sD_{T^\prime}) \\
 \end{array}\right],
\end{equation}
where $\Gamma$ is a contraction from $\sD_A$ into
$H^2(\sD_{T^\prime})$. Note that $\ga$ and $B$ in (\ref{defB})
define each other uniquely. Moreover, given (\ref{defB})  the second
identity in (\ref{rclt}) holds if and only if $\ga$ satisfies the
equation
\begin{equation}\label{fundeq}
E_{\sD_{T'}} D_{T^\prime}AR+S_{\sD_{T'}} \ga D_AR=\ga D_AQ.
\end{equation}
Therefore, with $U'$ as in (\ref{szns}), the RCL problem for
$\{A,T^\prime,U^\prime,R,Q\}$ is equivalent to the problem of
finding all contractions $\Gamma$ from $\sD_A$ into
$H^2(\sD_{T^\prime})$ such that (\ref{fundeq}) holds.

Equation (\ref{fundeq}) can be rewritten as an equation of the
form (\ref{cond1}). To see this one first observes that, because
of (\ref{intertw}), for each $h\in \sH_0$ we have
\begin{eqnarray*} \|D_A Q h\|^2 &=& \|Q h\|^2 - \|A Q h\|^2 \geq
\|R h\|^2 - \|T^\prime A R h\|^2\\
&=&  \|A R h\|^2 - \|T^\prime A R h\|^2 + \|R h\|^2 -\|A R h\|^2\\
&=& \|D_{T^\prime} A R h\|^2 + \|D_A R h\|^2.
\end{eqnarray*}
Hence the identity
 \[
\omega D_A Q h = \left[\begin{array}{c}
  D_{T^\prime} A R h \\
  D_A R h\\
\end{array}\right],\quad h \in \sH_0,
\]
uniquely defines a contraction $\omega$ from $\sF=\overline{D_A Q
\sH}$ into $\sD_{T^\prime} \oplus \sD_{A}$. We refer to this
contraction as the \emph{contraction underlying} the data set
$\{A,T^\prime,U^\prime,R,Q\}$. Using this contraction equation
(\ref{fundeq}) can equivalently be represented as
\begin{equation}\label{altfundeq}
E_{\sD_{T'}}\o_1+S_{\sD_{T'}}\Gamma\o_2=\Gamma|_{\sF},
\end{equation}
where $\omega_1$ is the contraction mapping $\sF$ into
$\sD_{T^\prime}$
 determined by the first component of $\omega$
and $\omega_2$ is the contraction mapping $\sF$ into $\sD_{A}$
determined by the second component of $\omega$. Summarizing the
above discussion we arrive at the following conclusion.

\medskip\noindent\textit{With $U'$ equal to the Sz.-Nagy-Sch\"{a}ffer isometric
lifting of $T'$, the RCL problem for $\{A,T^\prime,U^\prime,R,Q\}$
is equivalent to the problem of finding all contractions $\Gamma$
from $\sD_A$ into $H^2(\sD_{T^\prime})$ satisfying
$(\ref{altfundeq})$, where}
\[
\o=\left[\begin{array}{c}\o_1\\
\o_2
\end{array}\right]:\sF\to \left[\begin{array}{c}\sD_{T^\prime}\\
\sD_A
\end{array}\right]
\]
\textit{is the contraction underlying the given data set. Moreover,
the map $\ga \mapsto B$ determined by $(\ref{defB})$ provides a
one-to-one correspondence between the solutions of the interpolation
problem defined by $\o$ and the solutions of the RCL problem for
$\{A,T^\prime,U^\prime,R,Q\}$. In particular, any RCL problem
reduces to a problem of the type considered in the introduction.}

\section{Proofs of Theorems \ref{thmmain1} and \ref{thmmain2}}
\setcounter{equation}{0}

Throughout this section $\sU$ and $\sY$  are Hilbert spaces, $\sF$
is a subspace of $\sU$, and the operator $\o$ in (\ref{defom}) is
a contraction. We associate with $\o$ a lifting data set.

\begin{prop}
\label{omdata}Let $\o$ be a contraction as in $(\ref{defom})$. Put
\begin{eqnarray*}
&&\tilde{A}=\left[\begin{array}{cc}I_\sY&0\\0&0
\end{array}\right]\ \mbox{on} \ \left[\begin{array}{c}\sY\\\sU
\end{array}\right],  \quad \tilde{T}'=\left[\begin{array}{cc}0&0\\0&I_\sU
\end{array}\right]\ \mbox{on} \ \left[\begin{array}{c}\sY\\\sU
\end{array}\right]\\
\noalign{\vskip6pt} &&\tilde{R}=\left[\begin{array}{c}\o_1\\\o_2
\end{array}\right]:\sF\to\left[\begin{array}{c}\sY\\\sU
\end{array}\right],\hspace{.5cm}
\tilde{Q}=\left[\begin{array}{c}0\\\Pi_\sF^*
\end{array}\right]: \sF\to\left[\begin{array}{c}\sY\\\sU
\end{array}\right],\\
\noalign{\vskip6pt}
&&\tilde{U}'=\left[\begin{array}{ccc}0&0&0\\0&I_\sU&0\\E_\sY&0&S_\sY
\end{array}\right]\ \mbox{on} \
\left[\begin{array}{c}\sY\\\sU\\H^2(\sY)
\end{array}\right].
\end{eqnarray*}
Then $\{\tilde{A},\tilde{T}',\tilde{U'},\tilde{R},\tilde{Q}\}$ is a
data set, and the underlying contraction is precisely the given
contraction $\o$. Furthermore, $\tilde{U}'$ is the
Sz-Nagy-Sch\"affer isometric lifting of $\tilde{T}'$.
\end{prop}

Here $\Pi_\sF$ stands for the orthogonal projection of $\sU$ onto
$\sF$ viewed as a map from $\sU$ into $\sF$, and hence $\Pi_\sF^*$
is the canonical embedding of $\sF$ into $\sU$.

\medskip
\bpr The operators $\tilde{A}$ and $\tilde{T}'$ are orthogonal
projections and hence contractions. Observe that
$\tilde{T}'\tilde{A}$ and $\tilde{A}\tilde{Q}$ are both zero
operators. Furthermore, note that $R=\o$ is a contraction defined
on $\sF$ and $\tilde{Q}^*\tilde{Q}$ is the identity operator on
$\sF$. {}From these remarks we see that
\begin{equation}
\label{liftcond} \tilde{T}'\tilde{A}\tilde{R}=\tilde{A}\tilde{Q}
\quad \mbox{and}\quad \tilde{R}^*\tilde{R}\leq
\tilde{Q}^*\tilde{Q}.
\end{equation}
Next, observe that
\[
D_{\tilde{T}'}=\left[\begin{array}{cc}I_\sY&0\\0&0
\end{array}\right]\ \mbox{on} \ \left[\begin{array}{c}\sY\\\sU
\end{array}\right].
\]
Thus we can identify $\sD_{\tilde{T}'}$ with the space $\sY$. With
this identification in mind it is straightforward to check that
$\tilde{U}'$ is the Sz-Nagy-Sch\"affer isometric lifting of
$\tilde{T}'$. It follows that
$\{\tilde{A},\tilde{T}',\tilde{U}',\tilde{R},\tilde{Q}\}$ is a
data set. Notice that in this case the space $\sH_0$ appearing in
the definition of a data set is equal to the space $\sF$. Using
\[
D_{\tilde{A}}=\left[\begin{array}{cc}0&0\\0&I_\sU
\end{array}\right]\ \mbox{on} \ \left[\begin{array}{c}\sY\\\sU
\end{array}\right],
\]
we see that
\[D_{\tilde{A}}\tilde{Q}=\tilde{Q}, \quad
  D_{\tilde{T'}}\tilde{A}\tilde{R}=\o_1, \quad
  D_{\tilde{A}}\tilde{R}=\o_2.
\]
It is then easy to show that the contraction $\o$ in (\ref{defom})
is precisely the contraction underlying the data set
$\{\tilde{A},\tilde{T}',\tilde{U}',\tilde{R},\tilde{Q}\}$. \epr

\medskip
\noindent\textbf{Proof of Theorem \ref{thmmain1}.} Let $\o$ be a
contraction as in $(\ref{defom})$, and let
$\{\tilde{A},\tilde{T}',\tilde{U}',\tilde{R},\tilde{Q}\}$ be the
data set constructed in Proposition \ref{omdata}. Since $\o$ is the
contraction underlying this data set and $\tilde{U}'$ is
Sz-Nagy-Sch\"affer isometric lifting of $\tilde{T}'$, we know (see
the conclusion at the end of the previous section) that an operator
$\ga:\sU\to H^2(\sY)$ is a solution of interpolation problem defined
by the contraction $\o$ if and only if the operator
\begin{equation}
\label{rclom}
 B=\left[\begin{array}{c}\tilde{A}\\ \ga D_{\tilde{A}}
\end{array}\right]: \sY\oplus \sU\to \left[\begin{array}{c}\sY\oplus \sU\\H^2(\sY)
\end{array}\right]
\end{equation}
is a solution to the RCL problem for the data set
$\{\tilde{A},\tilde{T}',\tilde{U}',\tilde{R},\tilde{Q}\}$. Recall
(using canonical identifications) that $\sD_{\tilde{T}'}$ and
$\sD_{\tilde{A}}$ are equal to $\sY$ and $\sU$, respectively. But
then Theorem 1.1 in \cite{FtHK06b} tells use that $B$ in
(\ref{rclom}) is a solution to the RCL problem for the data set
$\{\tilde{A},\tilde{T}',\tilde{U}',\tilde{R},\tilde{Q}\}$ if and
only if the defining function $H$ of $\ga$ is given by
\[
H(\l)=\Pi_{\sY}Z(\l)(I-\l\Pi_{\sU}Z(\l))^{-1}, \quad \l\in\BD,
\]
where $Z$ is an arbitrary Schur  class function from
$\eS(\sU,\sY\oplus\sU)$ satisfying the constraint
$Z(\l)|_{\sF}=\o$ for each $\l\in\BD$. {}From these two ``if and
only if'' statements Theorem \ref{thmmain1} follows. \epr

\medskip
\noindent\textbf{Proof of Theorem \ref{thmmain2}.} Let $H\in
\eH^2(\sU,\sY)$ be a solution to the interpolation problem defined
by the contraction $\o$ in $(\ref{defom})$, and let $\ga $  from
$\sU$ into $H^2(\sY)$ be the operator defined by $H$. Then $\ga$ is
a contraction satisfying (\ref{cond1}). Hence for each $f\in \sF$ we
have
\begin{eqnarray}
\|D_\ga f\|^2&=& \|f\|^2- \|\ga f\|^2 =\|f\|^2- \|E_\sY\o_1
f\|^2-\|S_\sY\ga \o_2 f\|^2 \nonumber\\
&=&\|f\|^2- \|\o_1 f\|^2-\|\o_2 f\|^2+\|\o_2 f\|^2- \|\ga \o_2
f\|^2\label{omega12}\\
&\geq& \|D_\ga \o_2 f\|^2.\nonumber
\end{eqnarray}
The above calculation shows that there exists a (unique) contraction
$\om$ mapping $\sF_\ga=\ov{D_\ga\sF}$ into $\sD_\ga$ such that
\begin{equation}\label{defOG}
\om D_\ga|_\sF=D_\ga\o_2.
\end{equation}
Now let $\eS^\prime(\sD_\ga, \sD_\ga)$ be the set defined by
(\ref{parsetC1}). Let $C$ be a function in $\eS(\sD_\ga, \sD_\ga)$.
Using the identity (\ref{defOG}) we see that the function $C$
belongs to $ \eS^\prime(\sD_\ga, \sD_\ga)$ if and only if
$C(\l)D_\ga|_\sF=\om D_\ga|_\sF$. The latter identity can be
rewritten as $C(\l)|_{\sF_\ga}=\om$. We conclude that
\begin{equation}\label{parsetC2}
\eS^\prime (\sD_\ga, \sD_\ga)= \big\{ C\in \eS(\sD_\ga, \sD_\ga)\mid
C(\l)|_{\sF_\ga}=\om \ \mbox{for each  $\l\in \BD$}\big\}.
\end{equation}
 In other words, for the data set considered here the set (\ref{parsetC1})
 is precisely the set $\eS_\om(\sD_\ga, \sD_\ga)$ appearing in
 Theorem 1.2 of \cite{FtHK06b}. Next, recall that
\[
 B=\left[\begin{array}{c}\tilde{A}\\ \ga D_{\tilde{A}}
\end{array}\right]: \sY\oplus \sU\to \left[\begin{array}{c}\sY\oplus \sU\\H^2(\sY)
\end{array}\right]
\]
is a solution to the RCL problem for the data set
$\{\tilde{A},\tilde{T}',\tilde{U}',\tilde{R},\tilde{Q}\}$. But then,
using that the set (\ref{parsetC1}) is equal to the right hand side
of  (\ref{parsetC2}),  we can apply Theorem 1.2 in \cite{FtHK06b} to
complete the proof. \epr

\medskip

\begin{cor}
\label{coruni}Assume  that $\o$ in $(\ref{defom})$ is an isometry
such that $\o_2\sF$ is dense in $\sU$. Then the map $Z\mapsto H$
defined by Theorem \ref{thmmain1} is one-to-one, and
$(\ref{sols2})$ provides a proper parameterization of all
solutions to the interpolation problem defined by $\o$.
\end{cor}
\bpr Using the fact that $\o$ is an isometry, we see from
(\ref{omega12}) that the operator $\om$ in (\ref{parsetC2}) is
also an isometry. In particular, the space $\om\sF_\ga$ is a
closed subspace of $\sD_\ga$. Since $\o_2\sF$ is dense in $\sU$,
the space $D_\ga \o_2\sF$ is dense in $\sD_\ga$. By (\ref{defOG})
the space $D_\ga \o_2\sF$ is contained in $\om\sF_\ga$. Thus
$\om\sF_\ga$ is also dense in $\sD_\ga$. But $\om\sF_\ga$ is
closed in $\sD_\ga$, and therefore $\om_\ga\sF=\sD_\ga$. Hence
$\om$ is a unitary operator from $\sF_\ga$ onto $\sD_\ga$. This
implies that the set defined by the right hand side of
(\ref{parsetC2}), or equivalently the set (\ref{parsetC1}),
consists of one element only. Thus if $\o$ is an isometry and
$\o_2\sF$ is dense in $\sU$, then the map $Z\mapsto H$ defined by
Theorem \ref{thmmain1} is one-to-one. \epr

\medskip

{}From Theorem 1.2 in \cite{FtHK06b} it follows that the map
$C\mapsto Z_C$ given by (\ref{defZC}) and (\ref{defW}) is well
defined and induces a one-to-one map from the set (\ref{parsetC1})
onto the set of all $Z\in \eS(\sU,\sY\oplus\sU) $ satisfying
$Z(\l)|_{\sF}=\o$ for each $\l\in\BD$ and such that
$(\ref{sols2})$ holds.

\section{Proof of Theorem \ref{thmrepr}}
\setcounter{equation}{0} Throughout this section
$H\in\eH^2(\sU,\sY)$, and $\Theta\in \eS(\sE,\sU)$ is inner such
that $\Theta(0)=0$. Furthermore, $\sH=H^2(\sU)\ominus \Theta
H^2(\sE)$. Our aim is to prove Theorem \ref{thmrepr}. We begin with
some auxiliary results.

\begin{lem}\label{lemprel}
Let $\Phi=\l^{-1}\Theta$, and put $\sH_0=H^2(\sU)\ominus\Phi
H^2(\sE)$. Then
\begin{equation}
\label{formH}\sH=E_\sU\,\sU\oplus \l \sH_0, \quad
\sH=\sH_0\oplus\Phi E_\sE\,\sE.
\end{equation}
\end{lem}
\bpr As usual, given any inner function $\a\in \eS(\sX,\sY)$, we
shall denote the space $H^2(\sY)\ominus \a H^2(\sX)$  by $H(\a)$.
The two identities in (\ref{formH}), then follow from the rule
(see, e.g., Theorem X.1.9 in \cite{FF90}) that for two inner
functions $\a\in \eS(\sX,\sY)$ and $\b\in \eS(\sY,\sZ)$ we have
\[
H(\b\a)=H(\b)\oplus \b H(\a).
\]
Indeed, we apply this rule twice. First with $\a(\l)=\Psi(\l)$ and
$\b(\l)=\l I_\sU$, and next with $\a(\l)=\l I_\sE$ and $\b(\l)=
\Psi(\l)$. Note that in both cases $\b\a=\Theta$. With the first
choice of $\a$ and $\b$ we get the first identity in
(\ref{formH}), and the second choice yields the second identity in
(\ref{formH}). \epr

\medskip
The above lemma allows us to define the following auxiliary
operators:
\begin{equation}
\label{defRQ} R:\sH_0\to\sH,\quad Rh_0=h_0;\quad
Q:\sH_0\to\sH,\quad Qh_0=\l h_0\quad (h_0\in \sH_0).
\end{equation}
Note that the operators $R$ and $Q$ are isometries.

\begin{lem}\label{basiclem}
Let $\ga:\sH\to H^2(\sY)$ be a $($bounded linear$)$ operator, and
put
\begin{equation}
\label{defK}K(\l)=E_\sY^*(I-\l S_\sY^*)^{-1}\ga E_\sU, \quad \l\in
\BD.
\end{equation}
Then $K\in \eH^2(\sU,\sY)$, and the operator $\ga$  satisfies
$S_\sY\ga R=\ga Q$ if and only if for each $f\in \sH$ we have
\begin{equation}
\label{gaG} (\ga f)(\l)=K(\l)f(\l), \quad \l\in \BD.
\end{equation}
\end{lem}
\bpr According to the first identity in (\ref{formH}) the operator
$E_\sU$ maps $\sU$ into $\sH$. Thus ${\ga}E_\sU$ is well-defined,
and hence the same holds true for $K$.  Obviously, $K$ is an
$\sL(\sU,\sY)$-valued function which is analytic on $\BD$. Let
$K_n$ be the $n$-th coefficient of the Taylor expansion of $K$ at
zero. Take $u\in \sU$. Then $K_n
u=E_\sY^*(S_\sY^*)^n{\ga}E_\sU\,u=({\ga}E_\sU\,u)_n$, where
$({\ga}E_\sU\,u)_n$ is the $n$-th coefficient of the Taylor
expansion of the $\sY$-valued function ${\ga}E_\sU\,u$ at zero.
Therefore
\[
(\ga E_\sU u)(\l)=K(\l)u, \quad u\in\sU, \ \l\in \BD.
\]
Thus $K\in \eH^2(\sU,\sY)$, and $K$ is the defining function for
the operator $\ga E_\sU$.

Next assume that ${\ga}$ satisfies the intertwining relation
$S_\sY {\ga}R={\ga}Q$. {}From (\ref{defK}) we see that
$(\ref{gaG})$ holds for $f=E_\sU\,u$ with $u\in \sU$ arbitrary.
Indeed, for $f=E_\sU\,u$ we have $f(\l)\equiv u$, and hence
\[
K(\l)f(\l)=K(\l)u=E_\sY^*(I-\l
S_\sY^*)^{-1}({\ga}E_\sU\,u)=({\ga}E_\sU\,u)(\l)=({\ga}f)(\l).
\]
Using the first equality in (\ref{formH}), we see that it suffices
to prove $(\ref{gaG})$ for $f=\l h_0$ with $h_0\in \sH_0$. However
for such a function $f$ we have
\[
{\ga}f={\ga}Qh_0=S_\sY {\ga}Rh_0=S_\sY {\ga}h_0=\l {\ga}h_0,\quad
K(\l)f(\l)=\l K(\l)h_0(\l).
\]
Therefore, it suffices to prove $(\ref{gaG})$ for $h_0\in \sH_0$.

Take $h_0\in \sH_0$. Note that the identity operator on $H^2(\sU)$
is equal to $S_\sU S_\sU^*+P_\sU$, where $P_\sU$ is the orthogonal
projection of $H^2(\sU)$ onto $E_\sU\,\sU$. The first identity in
(\ref{formH}) shows that $E_\sU\,\sU\subset \sH$. Thus ${\ga}P_\sU$
is well defined. Since $\sH_0$ is invariant under $S_\sU^*$, we see
that $S_\sU S_\sU^*h_0$ belongs to $\sH$, and hence
\begin{eqnarray*}
{\ga}h_0&=&{\ga}S_\sU S_\sU^* h_0+{\ga}P_\sU h_0={\ga}QS_\sU^* h_0+{\ga}P_\sU h_0\\
&=&S_\sY {\ga}R S_\sU^* h_0+{\ga}P_\sU h_0=S_\sY {\ga} S_\sU^*
h_0+{\ga}P_\sU h_0.
\end{eqnarray*}
Since $h_0$ is an arbitrary element of $\sH_0$, we get
\[
{\ga}h_0 -S_\sY {\ga} S_\sU^* h_0  ={\ga}P_\sU h_0, \quad h_0\in
\sH_0.
\]
By induction, using $\sH_0$ is invariant under $S_\sU^*$\,, the
preceding identity  yields
\[
{\ga}h_0=\sum_{\nu=0}^n S_\sY^\nu ({\ga} P_\sU) (S_\sU^*)^\nu h_0 +
S_\sY^{(n +1)}{\ga} (S_\sU^*)^{(n+1)} h_0, \quad  h_0\in \sH_0,
\]
for $n=0,1,2, \ldots$. Now fix $\l \in \BD$. {}From (\ref{defK})
we know that
\[
({\ga} P_\sU f)(\l)=({\ga} E_\sU E_\sU^*\, f)(\l)=K(\l)E_\sU^*\,
f, \quad f\in H^2(\sU).
\]
Hence for $h_0\in \sH_0$ we have
\begin{eqnarray} &&({\ga}h_0)(\l)=\sum_{\nu=0}^n \l^\nu
K(\l)E_\sU^*\, (S_\sU^*)^\nu h_0+ \l^{n+1}\big({\ga}
(S_\sU^*)^{(n+1)} h_0\big)(\l)\nonumber\\
&&\quad =K(\l)\big(\sum_{\nu=0}^n \l^\nu E_\sU^*\, (S_\sU^*)^\nu
h_0\big) + \l^{n+1}\big({\ga} (S_\sU^*)^{(n+1)}
h_0\big)(\l).\label{partialsum}
\end{eqnarray}
 Note that for $n\to \iy$ the function ${\ga}
(S_\sU^*)^{(n+1)} h_0$ converges to zero in the norm of
$H^2(\sY)$, and hence the same holds true for $S_\sY^{(n +1)}{\ga}
(S_\sU^*)^{(n+1)} h_0$. It follows that the second term in the
right side of (\ref{partialsum}) converges to zero when $n\to
\iy$. Furthermore, for $n\to \iy$ the vector $\sum_{\nu=0}^n
\l^\nu E_\sU^*\, (S_\sU^*)^\nu h_0$ converges to $h_0(\l)$ in the
norm of $\sU$. Hence the first term in the right hand side of
(\ref{partialsum}) converges to $K(\l)h_0(\l)$ when $n\to \iy$.
Thus we have proved that (\ref{gaG}) holds.

To prove the converse implication. Assume that (\ref{gaG}) holds.
Let $h_0\in \sH_0$. Then for each $\l\in \BD$ we have
\begin{eqnarray*}
&&({\ga}Qh_0)(\l)=K(\l)\l h_0(\l)=\l K(\l)h_0(\l)=(S_\sY
{\ga}h_0)(\l)\\
\noalign{\vskip6pt} &&\hspace{2.1cm}=(S_\sY {\ga}Rh_0)(\l).
\end{eqnarray*}
Since $h_0$ is an arbitrary element in $\sH_0$, we see that $S_\sY
{\ga}R={\ga}Q$. \epr

\medskip
Next we put the problem into the setting of our alternative version
of the relaxed commutant lifting problem. Put $\sF=\l \sH_0$, and
define
\begin{equation}\label{defOm}
\o=\left[\begin{array}{c}\o_1\\ \o_2
\end{array}\right]:\sF\to \left[\begin{array}{c}\sY\\ \sH
\end{array}\right], \quad \o_1=0, \quad \o_2(\l h_0)=h_0 \
(h_0 \in \sH_0).
\end{equation}
Note that $\o_2Q=R$, and hence a contraction $\ga$ from $\sH$ into
$H^2(\sY)$ satisfies the intertwining relation $S_\sY\ga R=\ga Q$ if
and only if $S_\sY\ga\o_2=\ga|_\sF$. Since $\o_1=0$, we are now
ready to prove the main result.

\medskip\noindent\textbf{Proof of Theorem \ref{thmrepr}.}
Let $H\in \eH^2(\sU,\sY)$, and let us assume that the map
$f\mapsto Hf$ defines a contraction from $\sH$ into $H^2(\sY)$.
Denote this contraction by $\ga$. Then (\ref{defK}) holds with $H$
in place of $K$, and Lemma \ref{basiclem} shows that the
contraction $\ga$ satisfies the intertwining relation $S_\sY\ga
R=\ga Q$, and thus $S_\sY\ga\o_2=\ga|_\sF$. Since $\o_1=0$, we
know from Theorem \ref{thmmain1} that $H$ is given by
\begin{equation}\label{formH1}
H(\lambda)=F(\lambda)(I_\sH-\lambda G(\lambda))^{-1}E_\sU ,\quad
 \lambda\in\BD,
\end{equation}
where
\begin{equation}
\label{condsol} W=\left[\begin{array}{c}F\\ G
\end{array}\right]\in \eS(\sH, \left[\begin{array}{c}\sY\\ \sH
\end{array}\right]), \quad W(\l)|_\sF=\o\ \ (\l\in \BD).
\end{equation}

 Conversely, let $H\in \eH^2(\sU,\sY)$ be given by
(\ref{formH1}), where $F$ and $G$ are such that (\ref{condsol})
holds. Then, again using Theorem \ref{thmmain1}, we know that
there exists a contraction $\ga$ from $\sH$ into $H^2(\sY)$ such
that
\[
(\Gamma h)(\lambda)=F(\lambda)(I_\sH-\lambda
G(\lambda))^{-1}h,\quad h\in\sH,\ \lambda\in\BD.
\]
Moreover, $S_\sY\ga\o_2=\ga|_\sF$, and hence $\ga$ satisfies the
intertwining relation $S_\sY\ga R=\ga Q$.  It follows that
$H(\l)=(\ga E_\sU)(\l)$ for  each $\l \in \BD$. The fact that
$\ga$ satisfies the intertwining relation $S_\sY\ga R=\ga Q$
allows us to again apply  Lemma \ref{basiclem}. We conclude that
$(\ga h)(\l)=H(\l)h(\l)$ for each $h\in \sH$ and each $\l \in
\BD$. Thus the map $f\mapsto Hf$ induces a contraction from $\sH$
into $H^2(\sY)$ as desired.

{}From the previous results we see that it remains to show that
the representations given by (\ref{formH1}), (\ref{condsol}) and
by (\ref{functH}), with $Z\in \eS(\sU, \sY\oplus\sE)$, are
equivalent. Consider the spaces
\begin{equation}\label{sFsG}
\sF=\l \sH_0, \  \sG=\sH\ominus\sF=E_\sU\sU,\  \sF'=\sH_0, \
\sG'=\sH\ominus\sF'=\Phi E_\sE\sE.
\end{equation}

Let $F$ and $G$ be as in  (\ref{condsol}). Fix $z\in \BD$. Since
$G(z)|_\sF=\o_2$ and $\o_2$ is an isometry which maps $\sF$ onto
$\sF'$, we know that $G(z)\sG\subset\sG'=\Phi E_\sE\sE$. Thus
given $u\in \sU$, we have $G(z)E_\sU u=\Phi E_\sE e(z)$ for some
$e(z)\in \sE$.  Let $M_\Phi$ be the operator of multiplication by
$\Phi$ acting from $H^2(\sE)$ onto $H^2(\sU)$. The fact that
$\Phi$ is inner is equivalent to the statement that $M_\Phi$ in an
isometry. Put $C(z)=E_\sE ^*M_\Phi^*G(z) E_\sU$. Then
\[
C(z)u=E_\sE ^*M_\Phi^*G(z) E_\sU u=E_\sE ^*M_\Phi^*\Phi E_\sE
e(z)=E_\sE^* E_\sE e(z)=e(z).
\]
We conclude that $G(z)E_\sU=\Phi E_\sE C(z)$. {}From the
definition of $C(z)$ it is clear that $C(z)$ is a bounded linear
operator from $\sU$ into $\sE$. Moreover, $C(z)$ depends
analytically on $z\in \BD$.

{}From the result of the previous paragraph we know that
\[G(z)(E_\sU
u\oplus \l h_0)= h_0\oplus \Phi E_\sE C(z)= E_\sU v\oplus \l k_0,
\]
where
\[
v=E_\sU^*\Phi E_\sE C(z)u+ E_\sU^* h_0, \quad k_0=S_\sU^*\Phi
E_\sE C(z)u+ S_\sU^* h_0.
\]

Recall that $\sH=E_\sU\sU\oplus \l \sH_0$. Let $J$ be the operator
from $E_\sU\sU\oplus \l \sH_0$ to $\sU\oplus \sH_0$ defined by
$J(E_\sU u\oplus \l h_0)=u\oplus h_0$. Obviously, $J$ is unitary
and its inverse is given by $J^{-1}(u\oplus  h_0)=E_\sU u\oplus \l
h_0$.

It follows that relative to the direct sum decomposition
$\sU\oplus\sH_0$ the operator $JG(z)J^{-1}$ is given by the
following $2\times 2$ operator matrix:
\begin{equation}
\label{formG} JG(z)J^{-1}=\left[\begin{array}{cc}E_\sU^*\Phi E_\sE
C(z)&
E_\sU^*\\
\noalign{\vskip6pt} S_\sU^*\Phi E_\sE C(z)& S_\sU^*
\end{array}\right].
\end{equation}
But then
\[
J\big(I-zG(z)\big)J^{-1}=\left[\begin{array}{cc}I-z E_\sU^*\Phi
E_\sE C(z)& -z E_\sU^*\\
\noalign{\vskip6pt}-z S_\sU^*\Phi E_\sE C(z)& I-zS_\sU^*
\end{array}\right],
\]
and hence, using a Schur complement argument, we have
\[
J\big(I-zG(z)\big)^{-1}J^{-1}=\left[\begin{array}{cc}\Delta(z)^{-1}&
*\\
\noalign{\vskip6pt}*& *
\end{array}\right],
\]
where \begin{eqnarray*} \Delta(z)&=&I-z E_\sU^*\Phi E_\sE
C(z)-(-z E_\sU^*)(I-zS_\sU^*)^{-1}(-z S_\sU^*\Phi E_\sE C(z))\\
\noalign{\vskip6pt} &=& I-z E_\sU^*\Phi E_\sE C(z)+z
E_\sU^*(I-zS_\sU^*)^{-1}(I-z S_\sU^*-I)\Phi E_\sE C(z)\\
\noalign{\vskip6pt} &=& I-z E_\sU^*(I-zS_\sU^*)^{-1}\Phi E_\sE
C(z)\\
\noalign{\vskip6pt} &=& I-z \Phi(z)C(z)= I-\Theta(z)C(z).
\end{eqnarray*}
We also know that $F(z)|_\sF=0$. Thus $F(z)J^{-1}$ admits the
following representation
\begin{equation}
\label{formF} F(z)J^{-1}=\left[\begin{array} {cc}F_1(z)&0
\end{array}\right]:\left[\begin{array} {cc}\sU\\ \sH_0
\end{array}\right]\to \sY.
\end{equation}
But then we have
\begin{eqnarray*}
H(z)&=&F(z)\big(I-zG(z)\big)^{-1}E_\sU\\
\noalign{\vskip6pt}
&=&F(z)J^{-1}J\big(I-zG(z)\big)^{-1}J^{-1}JE_\sU\\
\noalign{\vskip6pt} &=&\left[\begin{array} {cc}F_1(z)&0
\end{array}\right]\left[\begin{array}{cc}\Delta(z)^{-1}&
*\\
\noalign{\vskip6pt}*& *
\end{array}\right]\left[\begin{array} {cc}I_\sU\\ 0
\end{array}\right]\\
&=&F_1(z)\Delta(z)^{-1}= F_1(z)\big(I-\Theta(z)C(z)\big)^{-1}.
\end{eqnarray*}

Let $\tau$ be  the canonical embedding of $\sH_0$ into $\sH$, that
is, $\tau$ is defined by $\tau h_0=h_0$. {}From (\ref{formG}) and
(\ref{formF}) it follows  that
\begin{equation}
\label{formFG}\left[\begin{array}{c}F(z)\\G(z)
\end{array}\right]J^{-1}=\left[\begin{array}{cc}F_1(z)&0\\
\noalign{\vskip6pt} \Phi E_\sE C(z)&\tau_{\sH_0}
\end{array}\right]:\left[\begin{array} {cc}\sU\\ \sH_0
\end{array}\right]\to \left[\begin{array} {cc}\sY\\ \sH
\end{array}\right].
\end{equation}
Since $\im \tau$ is perpendicular to $\Phi E_\sE\sE$ we see that
for $h=E_\sU u\oplus \l h_0$ we have
\[
\|F(z)h\|^2+\|G(z)h\|^2=\|F_1(z)u\|^2+\|\Phi E_\sE
C(z)u\|^2+\|h_0\|^2
\]
But $\|h\|^2=\|u\|^2+\|h_0\|^2$, and hence
\[
\|h\|^2-\big(\|F(z)h\|^2+\|G(z)h\|^2\big)=\|u\|^2-\big(\|F_1(z)u\|^2+\|\Phi
E_\sE C(z)u\|^2\big).
\]
Since multiplication by $\Phi$ and the map $E_\sE$ are isometries,
we conclude that
\[
W=\left[\begin{array}{c}F\\G
\end{array}\right]\in \eS(\sU,\left[\begin{array}{c}\sY\\ \sH
\end{array}\right]) \Longleftrightarrow Z=\left[\begin{array}{c}F_1\\C
\end{array}\right]\in \eS(\sU,\left[\begin{array}{c}\sY\\ \sE
\end{array}\right]).
\]
We have now shown that  the representations given by
(\ref{formH1}) and  (\ref{condsol}) imply those given by
(\ref{functH}), with $Z\in \eS(\sU, \sY\oplus\sE)$. The reverse
implication is obtained by reversing the arguments. \epr


\begin{thebibliography}{99}

\bibitem{ABL96}
D. Alpay, V. Bolotnikov, and Ph. Loubaton, On tangential $H_2$
interpolation with second order norm constraints, {\em Integral
Equations and Operator Theory} {\bf 24} (1996), 156--178.


\bibitem{ABL97}
D. Alpay, V. Bolotnikov, and Ph. Loubaton, On interpolation for
Hardy functions in a certain class of domains under moment type
constraints, {\em Houston Journal of Mathematics}, No. {\bf 3}
(1997), Vol. {\bf 23}, 539--571.


\bibitem{FF90}C. Foias and A. E. Frazho, \emph{The Commutant Lifting
Approach to Interpolation Problems,}  Operator Theory: Advances
and Applications, {\bf 44}, Birkh\"{a}user-Verlag, Basel, 1990.

\bibitem{FFK02a}
C. Foias, A.E. Frazho, and M.A. Kaashoek, Relaxation of metric
constrained interpolation and a new lifting theorem, {\em Integral
Equations and Operator Theory}, {\bf 42} (2002), 253--310.

\bibitem{FtHK06a}
A.E. Frazho, S. ter Horst, and M.A. Kaashoek, Coupling and relaxed
commutant lifting, {\em Integral Equations and Operator Theory},
{\bf 54} (2006), 33--67.

\bibitem{FtHK06b}
A.E. Frazho, S. ter Horst, and M.A. Kaashoek, All solutions to the
relaxed commutant lifting problem, {\em Acta Sci. Math (Szegged)}
{\bf 72} (2006), 299--318.

\bibitem{tH06a} S. ter Horst, Relaxed commutant lifting and a
relaxed Nehari problem: Redheffer state space formulas, to
submitted.

\bibitem{Sz.-NF} B. Sz.-Nagy and C. Foias, {\em Harmonic Analysis
of Operators on Hilbert Space,}  North Holland Publishing Co.,
Amsterdam-Budapest,  1970.

\bibitem{Sar89}D. Sarason, Exposed points in $H\sp 1$. I, in:
\emph{The Gohberg anniversary collection, Vol.~II}, OT
\textbf{41}, Birkh\"auser Verlag, Basel, 1989, pp. 485--496.



\end{thebibliography}
\end{document}